\documentclass{amsart}

\usepackage{amssymb}

\usepackage[all]{xy}

\newtheorem{thm}{Theorem}

\newtheorem{lem}[thm]{Lemma}

\newtheorem{prop}[thm]{Proposition}

\newtheorem{theorem}[thm]{Theorem}
\newtheorem{definition}[thm]{Definition}
\newtheorem{example}[thm]{Example}
\newtheorem{lemma}[thm]{Lemma}
\newtheorem{corollary}[thm]{Corollary}

\theoremstyle{definition}
\newtheorem{defn}[thm]{Definition}

\newtheorem{say}[thm]{}
\newtheorem{exmp}[thm]{Example}

\newtheorem{rem}[thm]{Remark}          

\newtheorem*{ack}{Acknowledgments}      

\newtheorem{defn-thm}[thm]{Definition--Theorem}  
\newtheorem{defn-lem}[thm]{Definition--Lemma}  

\theoremstyle{remark}

\renewcommand{\c}[0]{{\mathbb C}}  

\renewcommand{\o}[0]{{\mathcal O}} 
\newcommand{\z}[0]{{\mathbb Z}}

\renewcommand{\r}[0]{{\mathbb R}}

\newcommand{\p}[0]{{\mathbb P}}

\newcommand{\q}[0]{{\mathbb Q}}

\newcommand{\qtq}[1]{\quad\mbox{#1}\quad}

\newcommand{\pic}[0]{\operatorname{Pic}}
\newcommand{\gal}[0]{\operatorname{Gal}}

\newcommand{\inter}[0]{\operatorname{Int}}

\newcommand{\res}[0]{\operatorname{\mathcal R}} 
\newcommand{\dres}[0]{\operatorname{\mathcal D\mathcal R}} 
\newcommand{\link}[0]{\operatorname{link}}

\newcommand{\onto}[0]{\twoheadrightarrow}

\newcommand{\mor}[0]{\operatorname{Mor}} 
 
\newcommand{\tsum}[0]{\textstyle{\sum}}

\def\into{\DOTSB\lhook\joinrel\to}

\def\loccoh#1.#2.#3.#4.{H^{#1}_{#2}(#3,#4)}

\DeclareMathAlphabet{\mathchanc}{OT1}{pzc}%
                                {m}{it}

\def\CA{\mathcal A}
\def\CB{\mathcal B}
\def\CC{\mathcal C}
\def\CD{\mathcal D}
\def\CP{\mathcal P}

\def\CV{\mathcal V}
\def\CW{\mathcal W}

\def\C{\mathbb C}

\def\P{\mathbb P}
\def\Q{\mathbb Q}
\def\R{\mathbb R}
\def\Z{\mathbb Z}

\def\VV{\mathbf V}

\def\si{\sigma}

\def\<{\langle}
\def\>{\rangle}

\newcommand{\im}{\operatorname{im}}

\newcommand{\Span}{\operatorname{Span}}
\newcommand{\Faces}{\operatorname{Faces}}
\newcommand{\Resid}{\operatorname{Res}}
\newcommand{\Ob}{\operatorname{Ob}}
\newcommand{\Nerve}{\operatorname{Nerve}}

\begin{document}
\bibliographystyle{amsalpha}

\title[Fundamental groups of links]
{Fundamental groups of links of \\ isolated singularities}
\author{Michael Kapovich and J\'anos Koll\'ar}
\date{\today}       

\maketitle

Starting with   Grothendieck's proof of the local version of the 
Lefschetz hyperplane theorems \cite{SGA2}, it has been understood that 
there are strong parallels between the topology of smooth 
projective varieties and the
topology of links of isolated singularities. This relationship was 
formulated as one of the guiding principles in the monograph
 \cite[p.26]{gm-book}:
``Philosophically, any statement about the projective variety or its embedding
 really comes from a statement about the singularity at the point of the cone. 
Theorems about projective varieties should be consequences of more general 
theorems about singularities which are no longer required to be conical''.

The aim of this note is to prove the following, which we consider 
to be a strong exception to this principle.

\begin{thm} \label{link.thm}
For every finitely presented group $G$ there is an
isolated, 3-dimensi\-onal, complex singularity $\bigl(0\in X_{G}\bigr)$
with link $L_{G}$
such that  $\pi_1\bigl(L_{G}\bigr)\cong G$.
\end{thm}

By contrast, the fundamental groups of smooth projective varieties
are rather special; see \cite{abckt} for a  survey. 
 Even the fundamental groups of smooth  quasi projective
varieties are  quite restricted \cite{Morgan, kap-mil, Corlette-Simpson, DPS}.
This shows that germs of singularities can also be quite different
from quasi projective varieties.

We think of a  complex singularity $(0\in X)$ as a contractible Stein space 
sitting in some $\c^N$. Then its {\it link} 
is $\link (X):=X\cap {\mathbb S}^{2N-1}_{\epsilon}$, 
 an intersection of $X$ with a small $(2N-1)$-sphere
centered at $0\in X$. Thus $\link(X)$
 is a deformation retract of $X\setminus\{0\}$.

There are at least three natural ways to attach a fundamental group to an
isolated singularity $(0\in X)$. 
Let $p:Y\to X$ be a resolution of
the singularity with   simple normal crossing exceptional divisor $E\subset Y$. 
(That is, the irreducible components of $E$ are smooth and they intersects
transversally.)
We may assume
that $Y\setminus E\cong X\setminus\{0\}$. 
  The following 3 groups are
all independent of the resolution.
\begin{itemize}
\item $
\pi_1\bigl(\link(X)\bigr)=\pi_1\bigl(X\setminus\{0\}\bigr)=
\pi_1\bigl(Y\setminus E\bigr)$.
\item $\pi_1(Y)=\pi_1(E)$; we denote it by $\pi_1\bigl(\res (X)\bigr)$
to emphasize its independence of $Y$.
These groups were first studied  in  \cite{k-shaf, tak}. 
\item $\pi_1\bigl(D(E)\bigr)$ where $D(E)$ denotes the dual simplicial complex
 of $E$.
(That is, the vertices of $D(E)$ are the irreducible components  
$\{E_i:i\in I\}$ of  $E$ and for $J\subset I$ we attach a $|J|$-simplex for 
every irreducible component of $\cap_{i\in J} E_i$; see
 Definition \ref{DR.gen.defn} for details.)
We denote this group by $\pi_1\bigl(\dres (X)\bigr)$; it
is actually the main object of our interest.
\end{itemize}
There are natural surjections between these groups:
$$
\pi_1\bigl(\link(X)\bigr)=
\pi_1\bigl(Y\setminus E\bigr)\onto \pi_1(Y)=\pi_1\bigl(E\bigr)\onto 
\pi_1\bigl(D(E)\bigr).
$$
Usually neither of these maps is an isomorphism.
The kernel of $ \pi_1\bigl(Y\setminus E\bigr)\onto \pi_1(Y)$ is generated by
 loops around the irreducible components
of $E$, but the relations between these loops are not well understood;
see \cite{mumf} for some computations in the 2-dimensional case.

The kernel of $\pi_1(E)\onto \pi_1\bigl(D(E)\bigr) $ is generated by the
images of $\pi_1(E_i)\to \pi_1(E)$.  In all our examples the $E_i$ are 
simply connected,
thus $\pi_1(E)= \pi_1\bigl(D(E)\bigr) $. We do not investigate the difference
between these two groups in general.
\medskip

Our proof of Theorem \ref{link.thm} proceeds in two distinct steps.

Simpson showed in \cite[Theorem 12.1]{sim} that
 every finitely presented group $G$ is the
 fundamental group of a  singular  projective variety.
 He posed the question if this variety can be chosen to have 
 only simple normal crossing singularities. 
Our first result shows that this is indeed true:
 
\begin{thm}
For every  finitely-presented group $G$ there is a  complex, 
projective surface $S_G$ with simple normal crossing singularities only 
such that 
$\pi_1(S_G)\cong G$. 
\end{thm}
 
 Then we take the cone over $S_G$ to get an affine variety $C(S_{G})$
such that, essentially,   
$\pi_1\bigl(\res \bigl(C(S_{G})\bigr)\bigr)\cong G$.
(These cones have very non-isolated singularities
and therefore it is not clear that 
$\pi_1\bigl(\res \bigl(C(S_{G})\bigr)\bigr)$ really makes sense for them.)
Then we use the method of \cite{kol-exmp} to smooth
the singularities outside the origin while keeping the
fundamental group of the resolution unchanged.
(The smoothing yields an isolated singularity only in low dimensions;
 we get codimension 5 singularities in general.)
Thus we have a singularity $(0\in X_{G})$ such that
$\pi_1\bigl(\res \bigl(X_{G}\bigr)\bigr)\cong G$.
For some  choices of $S_G$ one can control the other two groups as well,
completing the proof of
 Theorem \ref{link.thm}.

\medskip

It is interesting to study the relationship between
the algebro-geometric properties of a singularity $(0\in X)$
and the  fundamental group of $\link(X)$. We prove the following
results for  rational
singularities in Section \ref{sec.ratsing}.

\begin{itemize}
\item  Let $(0\in X)$ be a rational singularity (\ref{rtl.defn}). Then
$\pi_1\bigl(\dres(X)\bigr)$ is  $\q$-super\-per\-fect,
that is,
$H_i\bigl( \pi_1\bigl(\dres(X)\bigr), \q\bigr)=0$ for $i=1,2$.
\item Conversely, for every finitely presented, $\q$-superperfect group $G$
there is a 6-dimensional rational singularity $(0\in X)$ 
such that 
$\pi_1\bigl(\dres(X)\bigr)=\pi_1\bigl(\res(X)\bigr)=\pi_1\bigl(\link (X)\bigr)
\cong G$.
\item Not every finite  group $G$ occurs as
$\pi_1\bigl(\dres(X)\bigr)$ for a 3-dimensional rational singularity.
(We do not know what happens in dimension 4 and 5.)\end{itemize}

\begin{say}[Open problems]\label{open.questions}
  Theorem  \ref{link.thm} and its proof 
raise many questions; here are just a few of them.

(\ref{open.questions}.1) Our examples show that links of isolated singularities
are  more complicated than {\em smooth} projective varieties. It would be interesting
to explore the difference in greater detail.

(\ref{open.questions}.2)  
Several steps of the proof use general position arguments and it is
probably impractical to follow it to get  concrete examples.
It would be nice to work out simpler versions in some key examples,
for instance  for Higman's group (as in Example \ref{higman})
and to understand the geometry of the resulting singularities completely.

(\ref{open.questions}.3) We have been focusing only on the fundamental group
of the dual simplicial complex  
$\pi_1\bigl(\dres(X)\bigr)$ associated to an isolated singularity $(0\in X)$. 
However, 
it is quite possible that for every finite simplicial complex ${\mathcal C}$ 
there is a
 an isolated singularity $(0\in X)$ such that $\dres(X)$ is 
homotopy equivalent to
 ${\mathcal C}$. Our results indicate that this may well be true but 
 our construction  yields non-isolated singularities
 and  only partial resolutions in general.

(\ref{open.questions}.4) 
Given an $n$-dimensional  manifold (possibly with boundary) $M$
our constructions  give a $(2n+1)$-dimensional singularity $(0\in X)$
such that $\dres(X)$ is homotopy-equivalent to $M$. It is reasonable to expect that there is an $(n+1)$-dimensional
singularity $(0\in Z)$ such that $\dres(Z)$ is homeomorphic to $M$.

(\ref{open.questions}.5)
For a complex algebraic variety $X$, its algebraic  fundamental group
$\pi^{\rm alg}_1(X)$ is the profinite completion of its
topological    fundamental group $\pi_1(X)$.
There are   examples where the natural map $\pi_1(X)\to \pi^{\rm alg}_1(X)$ is
not injective \cite{toledo, trento-2}, but in all such known cases
the image of  $\pi_1(X)\to \pi^{\rm alg}_1(X)$ is infinite and very large.
We now have examples of  isolated  rational singularities such that
$\pi_1(\link (X))$ is infinite yet $ \pi^{\rm alg}_1(\link (X))$ 
 is the trivial group, see Corollary \ref{cor:trivial}.

All the examples in  Theorem  \ref{link.thm}  can be realized
on varieties defined over $\q$. Thus they have an  
algebraic  fundamental group $\pi^{\rm alg}_1\bigl(\link (X_{\q})\bigr)$ which is an 
extension of
the above $\pi^{\rm alg}_1\bigl(\link (X)\bigr)$ and of the absolute Galois group
$\gal\bigl(\bar{\q}/\q\bigr)$. 
We did not investigate this extension, thus we
do not have a complete description of all possible groups
 $\pi^{\rm alg}_1\bigl(\link (X_{\q})\bigr)$.
\end{say}

\section{Polyhedral complexes}

A (convex) {\em Euclidean polyhedron} is a subset $P$ of $\R^n$ given by a finite collection of linear inequalities 
(some of which may be strict and some not). The polyhedron $P$ is {\em rational} if it can be given by linear inequalities with rational coefficients. The {\em dimension} of $P$ is its topological dimension, which is the same as the dimension of its affine span $\Span(P)$. (Recall that the empty set has dimension $-1$.) Note that we allow polyhedra which are unbounded and neither open nor closed. 
A {\em face} of $P$ is a subset of $P$ which is given by converting some of 
these non-strict inequalities to equalities. Define the set $\Faces(P)$ to be the set of faces of $P$. 
The {\em interior} $\inter(P)$ of $P$ is the topological interior of $P$ in $\Span(P)$. Again, 
$\inter(P)$ is a  Euclidean polyhedron. We will refer to $\inter(P)$ as an {\em open polyhedron}.

An (isometric)  {\em morphism} of two polyhedra is an isometric map $f: P\to Q$ so that 
$f(P)$ is  a face of $Q$. A morphism is {\em rational} if it is a restriction of a rational affine map. 


\begin{definition}
A {\em (Euclidean) polyhedral complex}  is a small category $\CC$ whose objects are convex polyhedra and 
morphisms are their isometric morphisms satisfying the following axioms:

1)  For every $c_1\in \Ob(\CC)$ and every face $c_2$ of $c_1$, $c_2\in \Ob(\CC)$,  the inclusion map $\iota: c_1\to c_2$ is a morphism of $\CC$.  

2) For every $c_1, c_2\in \Ob(\CC)$ there exists at most one morphism $f=f_{c_2,c_1}\in \mor(\CC)$ so that $f(c_1)\subset c_2$. 

A polyhedral complex is {\em rational} if all its objects and morphisms are rational.  

\end{definition}

Analogously, one defines spherical, hyperbolic, affine, projective, etc., polyhedral complexes, but we will not need these concepts. 
Thus, a {\em polyhedral complex} for us will always mean a Euclidean polyhedral complex.  

Given a finite rational polyhedral complex $\CC$, by scaling one obtains an {\em integral} polyhedral cell complex $\CC'$, where the 
polyhedra and morphisms are integral. 

\begin{example}
Every simplicial complex $Z$ corresponds canonically to a Euclidean polyhedral complex ${\mathcal Z}$: Identity each $k$-simplex in $Z$ with the standard Euclidean simplex  in $\R^{k+1}$. 
\end{example}

Objects of a polyhedral complex $\CC$ are called {\em faces} of $\CC$ and the morphisms of $\CC$ are called {\em incidence maps} of $\CC$. 
A {\em facet} of $\CC$  is a face $P$ of $\CC$ so that 
for every morphism $f: P\to Q$ in $\CC$, $f(P)=Q$.  A {\em vertex} of $\CC$ is a zero-dimensional face. The {\em dimension} 
$\dim(\CC)$ of $\CC$ is the supremum of dimensions of faces of $\CC$. 
A polyhedral complex $\CC$ is called {\em pure} if the dimension function is constant on the set of facets of $\CC$; 
the constant value in this case is the dimension of $\CC$. A {\em subcomplex} of  $\CC$  is a full subcategory of $\CC$.   
If $c$ is a face of a complex $\CC$ then $\Resid_{\CC}(c)$, the {\em residue} of $c$ in $\CC$, is the minimal 
subcomplex of $\CC$ containing all faces  $c'$ such that there exists an incidence map $c\to c'$. For instance, if $c$ is a vertex of 
$\CC$ then its residue is the same as the {\em star} of $c$ in $\CC$; however, in general these are different concepts. 
 
 \medskip 
We generate the equivalence relation $\sim$ on a polyhedral complex $\CC$ by declaring that $c\sim f(c)$, where $c\in \Ob(\CC)$ and 
$f\in \mor(\CC)$. This equivalence relation also induces the equivalence relation $\sim$ on points of faces of $\CC$.

If $\CC$ is a polyhedral complex, its {\em poset} $\operatorname{Pos}(\CC)$ is the partially ordered set $\Ob(\CC)$ with the relation $c_1\le c_2$ 
iff $c_1\sim c_0$ so that $\exists f\in \mor(\CC), f: c_0\to c_2$. 

We define the {\em underlying space} or {\em amalgamation} $|\CC|$ of a polyhedral complex $\CC$ as the topological space which is obtained from the disjoint union 
$$
\amalg_{c\in \Ob(\CC)}\ c
$$
by identifying points using the equivalence relation: $\sim$. We equip  $|\CC|$ with the quotient topology.




\begin{definition}\label{difference}
If $\CC$ is a polyhedral complex and $\CB$ is its subcomplex. For $c\in \Ob(\CC)$ define the polyhedron 
$$
 c':=c\setminus \bigcup_{b\le c, b\in \CB}  f(b), \hbox{~~where~~}  f: b\to c, f\in \mor(\CC). 
$$
For a morphism $f\in \mor(\CC)$, $f: c_1\to c_2$, we set $f': c_1'\to c_2'$ be the restriction of $f$. 
We define the {\em difference complex} $\CC - \CB$ as the following polyhedral complex:
$$
\Ob(\CC - \CB)=\{ c': c \in \Ob(\CC)\},   
$$
$$
\mor(\CC - \CB)= \{f': c_1'\to c_2',  \hbox{where~} f\in \mor(\CC), f: c_1\to c_2\}. 
$$
\end{definition}

Note that $\CC - \CB$ need not be a subcomplex of $\CC$. 

\medskip 
We will be exclusively interested in finite Euclidean polyhedral complexes (i.e., complexes of finite cardinality), 
where the underlying space $|\CC|$ is connected.  

\begin{example}\label{E1}
For us the most important  Euclidean polyhedral complexes are  obtained by subdividing (a domain in) 
$\R^N$ into convex polyhedra. Let $U\subset \R^N$ be an open subset and $\CV$ be a partition of $U$ in convex Euclidean polyhedra so that: 

1)  For every $c_1\in \CV$ and every face $c_2$ of $c_1$, $c_2\in \CV$. 

2) For every two polyhedra $c_1, c_2\in \CV$, $c_1\cap c_2\in \CV$.   

Then $\CV$ becomes the set of faces of a polyhedral complex (again denoted by $\CV$ by abusing the notation), 
where inclusions of faces are the incidence maps.   
\end{example}

\begin{example}
Let $\Delta^m$ be the closed Euclidean $m$-simplex. The (simplicial) cell complex of faces of $\Delta^m$ will be denoted $\CC(\Delta^m)$. 
\end{example}

\begin{definition}
Let $\CC$ be a pure $n$-dimensional polyhedral complex. 
The {\em nerve} $\Nerve(\CC)$ of $\CC$ is the simplicial complex whose vertices are facets of $\CC$ 
(the notation is $v=c^*$, where $c$ is a facet of $\CC$); 
distinct vertices $v_0=c_0^*,...,v_k=c_k^*$ or $\Nerve(\CC)$ span a $k$-simplex if there exists an $n-k$-face $c$ of $\CC$ and 
incidence maps $c\to c_i, i=0,...,k$. The simplex $\sigma=[v_0,...,v_k]$ then is said to be {\em dual} to the face $c$.    
\end{definition}

\begin{lemma}
If $\CC$ is finite then $|\CC|$ is homotopy-equivalent to $|\Nerve(\CC)|$. 
\end{lemma}
\proof Notice that the image of each face of $\CC$ in $|\CC|$ is contractible. Therefore, one can thicken each facet $c$ of $\CC$ to 
an open contractible subset $U(c)\subset |\CC|$ so that: 

(1) The collection of open sets $U(c)$, $c$'s are facets of $\CC$, is a covering of $|\CC|$. 

(2) For every $k+1$-tuple of facets $c_0,...,c_k$, the intersection 
$$
U(c_0)\cap ... \cap U(c_k) 
$$
is nonempty iff $[c_0^*,...,c_k^*]$ is a simplex in $\Nerve(\CC)$. 

(3) Each intersection $U(c_0)\cap ... \cap U(c_k)$ as above is contractible. 

Now, the assertion becomes the standard fact of algebraic topology, see e.g. \cite{Hatcher}. \qed 

\medskip 
Note that, in general, a face can have more than one, or none, dual simplex, except that each facet is dual to the unique vertex.  

\begin{definition}
A polyhedral complex $\CC$ is {\em simple} if: 

(1) $\CC$ is pure, $\dim(\CC)=n$,  

(2) For $k=0,...,n$ and every $k$-face $c$ of $\CC$, 
$\Nerve(\Resid_{\CC}(c))$ is isomorphic to the complex $\CC(\Delta^{n-k})$. 


\end{definition}

It is easy to see that each face $c$ of a simple $n$-dimensional 
complex $\CC$ is dual to a unique simplex $c^*$ in $\Nerve(\CC)$. Moreover, $\dim(c)+\dim(c^*)=n$.  

\begin{lemma}
If $\CA$ is a simple polyhedral complex and $\CB$ is its subcomplex, then the complex $\CC:=\CA - \CB$ is again simple. 
\end{lemma} 
\proof It is easy to see that $\CC$ is is pure of the same dimension as $\CA$ and 
for each face $c$ of $\CC$, the poset of $\Resid_{\CC}(c)$ is isomorphic to the poset $\Resid_{\CA}(a)$, where $c=a'$ (see Definition \ref{difference}). 
\qed


\section{Voronoi complexes in $\R^N$}

In what follows, we will use the notation $d$ for the 
Euclidean metric on $\R^N$.

\begin{defn}
Let $Y\subset \R^N$ be a finite subset. 
The {\em Voronoi tessellation} 
$\CV(Y)$ of $\R^N$ associated with $Y$ is defined by:  For each $y\in Y$ take the {\em Voronoi cell} 
$$V(y):=\{x\in \R^N: d(x,y)\le d(x,y'), \forall y'\in Y\}.$$ 
Thus, each cell $V(y)$ is given by the collection of non-strict linear inequalities $d(x,y)\le d(x,y')$, i.e., 
$$
2 (y'-y)\cdot x \le y'\cdot y' - y\cdot y. 
$$ 
Then each cell $V(y)$ is a closed (possibly unbounded) polyhedron in $\R^N$. 
Every $V(y)$ is rational provided that $Y\subset \Q^N$. The union of Voronoi cells is the entire $\R^N$. 
We thus obtain the polyhedral complex, called the {\em Voronoi complex}, $\CV(Y)$ using the faces $V(y)$ as in Example \ref{E1}. 
\end{defn}

Not every Voronoi complex is simple, but most of them are.
 In order to make this precise, we  consider  ordered finite subsets of $\R^N$; thus, every $k$-element subset 
becomes a point in $\R^{kN}$. 

\begin{lemma}\label{simple}
 $k$-element subsets  $Y\subset \R^N$ (resp.\ $Y\subset \q^N$)
whose Voronoi complex $\CV(Y)$ is simple
are open and dense in $\R^{kN}$  (resp.\ $\q^{kN}$).
\end{lemma}
\proof

For a  subset $Y\subset \R^N$, failure of simplicity of ${\CV}(Y)$ means that there exists an $m$-element subset $W\subset Y$ so that the set of affine hyperplanes  
$$
H_{y_i,y_j}=\{x| d(y_i,x)=d(y_j,x)\}, y_i, y_j\in W$$ 
has non-transversal intersection in $\R^N$. The subset $\Sigma_{m,N}\subset \R^{mN}$ of such $W$'s is closed and has empty interior. The 
Lemma follows from density of rational $k$-element subsets. \qed

\medskip 
{\bf Delaunay triangulations.} Dually, one defines the {\em Delaunay simplicial complex} $$\CD(Y)=\Nerve(\CV(Y)),$$ i.e., vertices of this complex are points 
of $Y$, vertices $y_0,...,v_k$ span a $k$-simplex in $\CD(Y)$ iff 
$\cap_{i=0}^k D(y_i) \ne \emptyset$.
We have the canonical affine map $\eta: \CD(Y)\to \R^N$ which is the identify on $Y$. 

The proof of the following theorem can be found in \cite[Thm.2.1]{Fortune}:

\begin{theorem}\label{Delaunay}
1. If $\CV(Y)$ is simple then the map $\eta: |\CD(Y)|\to \R^N$ injective.

2. The image of the latter map is the closed convex hull 
$\operatorname{Hull}(Y)$ of the set $Y$. 

\end{theorem}
 
\noindent The Euclidean simplicial complex $\eta(\CD(Y))$ is called the {\em Delaunay triangulation} of  $\operatorname{Hull}(Y)$  
with the vertex set $Y$.


\medskip
{\bf Voronoi complexes associated with smooth submanifolds in $\R^N$.} 
Let $M$ be a subset of $\R^N$ and $\epsilon>0$. A set $Y\subset \R^N$ is said to be $\epsilon$--dense in $M$ if 
every point $x\in M$ is within distance $<\epsilon$ from a point of $Y$. (Note that $Y$ need not be contained in $M$.) 
By compactness and Lemma \ref{simple}, 
for every bounded subset of $\R^N$, $\epsilon>0$, there exists a finite simple rational subset $Y\subset \R^N$ which is $\epsilon$-dense in $M$.

\begin{theorem}\label{cairns}
(Cairns \cite{Cairns}.) 
Let $M$ be a $C^2$-smooth closed submanifold  of $\R^N$. Then there exists $\epsilon_M>0$ so that for every $\epsilon\in (0, \epsilon_M]$ 
the following holds:

Let $Y$ be a finite subset of $\R^N$ which is $\epsilon$-dense in $M$. For every face $c\in \CV(Y)$ define the set $c_M:=c\cap M$.

Then each $c_M$ is a (topological) cell in $M$ and the collection of cells $c_M, c\in \CV(Y)$, 
is a cellulation of $M$. 
\end{theorem}

Note that for a generic choice of $Y$ the intersections $c_M$ are transversal and, hence, 
$$
\dim(c_M)= \dim(c) + \dim(M) -N . 
$$

Examination of Cairns's proof of this theorem shows that it can be repeated verbatim to prove the following

\begin{theorem}
Let $S$ be a compact codimension $0$ submanifold  of $\R^N$ with $C^2$-smooth boundary $M$. 
Then for $\epsilon_M>0$ as above, and for every $\epsilon\in (0, \epsilon_M]$ 
the following holds:

Let $Y$ be a finite subset of $\R^N$ which is $\epsilon$-dense in $S$. For every face $c\in \CV(Y)$ define the set $c_S:=c\cap S$. 
Then each $c_S$ is a (topological) cell in $S$ and the collection of cells $c_S, c\in \CV(Y)$, 
is a regular cellulation $\CV_S$ of $S$. 

Again, for a generic choice of $Y$, $\dim(c_M)= \dim(c)$ unless $c_M=\emptyset$. 
\end{theorem}

We observe that $\{c\in \CV(Y): c_S=\emptyset\}$ is the face set of a subcomplex $\CW(Y)$ of $\CV(Y)$. 
We then let $\CC_S:= \CV(Y) - \CW(Y)$.

Clearly, $\CC_S$ is pure, $N$-dimensional and $\Nerve(\CC_S)$ is isomorphic to $\Nerve(\CV_S)$. 
Therefore, $|\Nerve(\CC_S)|$ is homotopy-equivalent to $|\Nerve(\CV_S)|\cong S$. 


By Lemma \ref{simple}, we can assume that $Y$ is simple, rational and generic; 
then $\CC_S$ is again simple, rational and $|\CC_S|$ is still is homotopy-equivalent to $S$. 
 
 \begin{corollary}\label{normal}
 Given $S\subset \R^N$ as above, there exists a simple rational polyhedral complex $\CC=\CC_S$ with simple parasites
so that $|\CC|$ is homotopy-equivalent to $S$. 
 \end{corollary}

Faces of $\CC_S$, in general, are not closed polyhedra.

\section{Euclidean thickening of simplicial complexes}\label{thickening}
  
 We are grateful to 
 Frank Quinn for 
 leading us to the following reference: 
  
 \begin{theorem}
(Hirsch \cite{Hirsch}) Let $Z$ be a finite simplicial complex in a smooth manifold $X$.  
Then there exists a codimension 0 compact submanifold $S\subset X$ with smooth boundary which is homotopy-equivalent to $|Z|$. 
\end{theorem}

\begin{corollary}\label{walker}
For every $n$-dimensional finite simplicial complex $Z$ there exists a codimension 0 compact submanifold $M\subset \R^{2n+1}$ with smooth boundary, homotopy-equivalent to $|Z|$.
\end{corollary}
\proof Below is a simple and self-contained proof of this corollary communicated to us by Kevin Walker via Mike Freedman. We will think of  $Z$ as a cell complex and will construct $M$ by induction on skeleta of $Z$. Let $B_v$ denote pairwise disjoint intervals in $\R$,  
where $v\in Z^{(0)}$ and set
$$
S_0:= \amalg_{v} B_v. 
$$
Then the map $Z^{(0)}\to S_0$ sending each $v$ to $B_v$ is a homotopy-equivalence. Suppose that we constructed a codimension $0$ 
submanifold $M_k\subset \R^{2k+1}$ with smooth boundary and a homotopy-equivalence $h: Z^{(k)}\to M_k$. Let $f_i$  denote the 
the attaching maps $\partial B^{k+1}\to Z$ of $(k+1)$-cells of $Z$. The maps $f_i$ define elements of $\pi_{k}(M_k)$. 
Since the dimension of the manifold $M_k$ is $2k+1$, 
these homotopy classes can be realized by  smoothly embedded $k$-sphere $s_i$ with trivial normal bundle in $\R^{2k+1}$, see 
\cite[\S 6]{Kervaire-Milnor}. A priori, the spheres $s_i$ may not be even homotopic to 
spheres contained in  the boundary of $M_{k}$. However, we replace $M_k$ with $M_k'\subset \R^{2k+3}$, 
obtained from $M_k\times B^2\subset \R^{2k+3}$ by ``smoothing the corners.'' Then $s_i$ can be chosen in the boundary of $M_k'$ and, 
for the dimension reasons, it bounds a smoothly embedded $(k+1)$-disk in $\R^{2k+3}\setminus M_k'$. Then we attach the handle 
$H_i\cong D^{k+1}\times D^{k+2}\subset \R^{2k+3}$ to $M_k'$ along $s_i$, so that this handle 
intersects $M_k'$ only along a tubular neighborhood of $s_i$ in $\partial M'$. Moreover, we can assume that distinct handles   $H_i$ 
are pairwise disjoint. We again smooth the corners after the handles are attached. 
Let $M_{k+1}\subset \R^{2k+1}$ be the codimension $0$ submanifold with smooth boundary resulting from attaching these handles and 
smoothing the corners. Then, clearly, the homotopy-equivalence $Z^{(k)}\to M_k$ extends to a 
homotopy-equivalence $Z^{(k+1)}\to M_{k+1}$. \qed

\medskip 
By combining Hirsch's theorem with Corollary \ref{normal} we obtain: 

\begin{corollary}\label{Hirsh+Cairns}
Given a finite $n$-dimensional simplicial complex $Z$, there exists a finite  
 simple $(2n+1)$-dimensional rational Euclidean polyhedral complex $\CC$ so that $|\CC|$ is 
homotopy-equivalent to $|Z|$. 
\end{corollary}

We note that the dimension of $\CC$ in this corollary can be easily reduced to $N=2n$. For instance, 
by a theorem of Stallings (see \cite{DR})  
$Z$ is homotopy-equivalent to a finite simplicial complex $W$ which embeds in $\R^{2n}$. One can improve on this estimate even further 
as follows. Suppose that $Z$ is a finite simplicial complex which admits an {\em immersion} $j: |Z|\to \R^N$. Then taking pull-back of an 
open regular neighborhood of $j(|Z|)$ via $j$ one obtains an open smooth locally-Euclidean $N$-dimensional manifold $X$ which is 
homotopy-equivalent to $|Z|$. 

\begin{definition}
If $Z$ is a simplicial complex, then a locally-Euclidean Riemannian manifold $X$ is called a {\em Euclidean thickening} of $Z$ if 
there exists an embedding $|Z|\to X$ which is a homotopy-equivalence. We say that  $X$ is {\em rational}, resp. {\em integral} if there exists a smooth atlas on $X$ with transition maps that belong to $\Q^n \rtimes GL(n,\Q)$, resp. $\Z^n \rtimes GL(n,\Z)$, where 
$n=dim(X)$. 
\end{definition}

Note that if $X$ is an $n$-dimensional  locally-Euclidean manifold which admits an isometric immersion in $\R^n$, then $X$ is integral. 

\medskip 
Suppose that $X$ is a Euclidean thickening of $Z$. 
Applying Hirsch's   theorem above to the embedding $|Z|\subset X$, we obtain an $N$-dimensional 
smooth manifold with boundary $S\subset X$ homotopy-equivalent to $|Z|$. 
Even though such $X$ is not isometric to $\R^N$ (it is typically incomplete and, moreover, need not embed in $\R^N$ isometrically), the arguments in the proof of Cairns's theorem \ref{cairns} are local 
and go through if we replace $\R^N$ with $X$. 
We thus obtain

\begin{corollary}
Suppose that $Z$ is a finite simplicial complex and $X$ is an $N$-dimensional Euclidean thickening of $Z$. Then there exists a simple 
$N$-dimensional Voronoi complex $\CC$ 
so that $|\CC|$ is homotopy-equivalent to $Z$. Moreover, if $X$ is rational, resp. integral, 
the complex $\CC$ can be taken rational, resp. integral. 
\end{corollary}

\section{Complexification of Euclidean polyhedral complexes}
\label{sec.4}


\begin{defn} Let $\VV$ denote either the category of varieties
(over a fixed field $k$) or the category of topological spaces.

Let $\CC$ be a finite polyhedral complex. 
A {\it  $\VV$-complex based on $\CC$} is a functor $\Phi$ 
from $\CC$ to $\VV$ so that 
 morphisms $c_i\to c_j$ go to closed embeddings 
$\phi_{ij}:\Phi(c_i)\into \Phi(c_j)$. By abuse of terminology, we will sometimes refer to the image category $im(\Phi)$ as a $\VV$-complex based on $\CC$. 

We call the functor $\Phi$ {\it strictly faithful} if the following holds:

 If $x_i\in \Phi(c_i)$, $x_j\in \Phi(c_j)$ and $\phi_{ik}(x_i)=\phi_{jk}(x_j)$ 
for some $k$ then there is an $\ell$
and $x_{\ell}\in \Phi(c_{\ell})$ such that
$\phi_{\ell i}(x_{\ell})=x_i$ and $\phi_{\ell j}(x_{\ell})=x_j$.

The relation $x_i\sim \phi_{ij}(x_i)$ for every $i, j$ and $ x_i\in X_i$
generates an equivalence relation
on the points of $\amalg_{i\in I} \Phi(c_i)$, also denoted by $\sim$. 
\end{defn}

In the category of topological spaces, the direct limit $\lim\Phi(\CC)$
of the diagram $\Phi(\CC)$  exists and its points are identified with
 $\bigl(\amalg_{i\in I} \Phi(c_i)\bigr)/\sim $.

For example, suppose that $\Phi_{taut}$ is the tautological functor which identifies each face of $\CC$ with the corresponding underlying topological space. 
Then $\lim \Phi_{taut}(\CC)$ is nothing but $|\CC|$. 

In general, Proposition 3.1 in \cite{Corson} proves the following.

\begin{lemma}\label{pi_1-cong}
Suppose that $\Phi$ is strictly faithful and $\Phi(\CC)$ consists of cell complexes and cellular maps of such complexes. Then 
 $\pi_1\bigl(\lim\Phi(\CC)\bigr)\cong \pi_1(|\CC|)$ provided that each $\Phi(c), c\in \Ob(\CC)$ is $1$-connected. \qed
\end{lemma}

In the category of varieties direct limits usually do not exist;
we deal with this question in Section \ref{sec.5}. 
Thus for now assume that $\Phi(\CC)$ has a direct limit $\lim\Phi(\CC)$
in  the category of varieties.
There is a natural surjection 
$$
\bigl(\amalg_{i\in I} \Phi(c_i)\bigr)/\sim \quad \to \quad \lim\Phi(\CC).
$$
If this map is a bijection, we say that $\lim\Phi(\CC)$ is
an {\it algebraic realization} of $|\CC|$.

\medskip Our next goal is to describe two constructions of complexes of varieties based on polyhedral complexes. 

\medskip
{\bf Projectivization.} Suppose that $\CC$ is a Euclidean polyhedral complex.

We first construct a {\em projectivization} $\CP=\CP(\CC)$ of the Euclidean polyhedral complex ${\CC}$. We regard each face $c$ of $\CC$ 
as a polyhedron in $\R^N$.
 The complex affine span, $\Span(c)$, of $c$ is a linear subspace of $\C^N$. Let $P_c$ denote its
 projective completion in $\P^N$.  
Note that, technically speaking, different faces $c$ can yield the same space $P_c$ if their affine spans are the same. To avoid these issues, we set $\CP_c:=P_c\times \{c\}$. For every morphism  of $\CC$, $f_{c_2,c_1}: c_1\to  c_2$, we have a unique linear embedding $F_{c_2,c_1}: \CP_{c_1}\to \CP_{c_2}$.
 Thus, we obtain the functor
$$
\CP: \CC \to ~~\hbox{Varieties}
$$
which sends each $c\in \Ob(\CC)$ to $\CP_c$ and each morphism $f_{c_2,c_1}$ to $F_{c_2,c_1}$. 

We refer to $\CP$  
as a  {\em complex of projective spaces} based on the complex ${\CC}$. 

Observe, however, that $im(\CP)$ does not (in general)   
accurately capture the combinatorics of the complex ${\mathcal C}$, i.e., the corresponding functor $\Psi_{var}=\CP$ need not be 
strictly faithful. For instance, in complex projective space any two
hyperplanes intersect but a  polyhedral complex usually has
disjoint codimension 1 faces. More generally, any intersection
$\cap_i P_{c_i}$ that is not  equal to the 
 projective span $P_c$ of some face $c\in \Ob(\CC)$
shows that $\Psi_{var}=\CP$ is not
strictly faithful.

\begin{defn}[Parasitic intersections]\label{paras.int.defn}
Let  $\sigma:=(c_1, c_2,...,c_k)$ be a tuple of faces
incident to a face $c$. Consider the intersections 
$$
I_{c, \sigma}:= \cap_{i=1}^k F_{c,c_i}(\CP_{c_i})\subset \CP_{c}
$$
such that there is no face $c_0$ such that
$I_{c, \sigma}=F_{c,c_0}(\CP_{c_0})$ and $c_0$ is incident to all the 
$c_1, c_2,...,c_k$. Then the subspace $I_{c, \sigma}\subset \CP_{c}$ 
is called a {\em parasitic intersection} in $\CP_{c}$.

Note, however, that this collection of parasitic intersections in spaces $\CP_c, c\in Ob(\CC)$, 
is not stable under applying morphisms $F_{c',c}$ and taking preimages under these morphisms. 
We thus have to {\em saturate} the collection of parasitic intersections using the morphisms $F_{c,c'}$. 
This is done as follows. Let $T$ denote the pushout of the category $im(\CP_{top})$, where we regard each $\CP_b, b\in Ob(\CC)$,   
as a topological space, so the pushout exists. Then for each $a\in Ob(\CC)$ we have the (injective) projection map 
$\rho_a: \CP_a\to T$.  For each parasitic intersection $I_{c,\si}\subset \CP_c$, we define
$$
I_{c,\si,a}:=\rho_a^{-1} \rho_c(I_{c,\si}). 
$$
We call such $I_{c,\si,a}$ a {\em parasitic subspace} in $\CP_a$. It is immediate that each parasitic subspace in $\CP_a$ 
is a projective space linearly embedded in $\CP_a$. With this definition, the collection of parasitic subspaces $I_{c,\si,a}$ 
is stable under taking images and preimages of the morphisms $F_{c,c'}$. 
\end{defn}

Note that the maps $\rho_c$ induce an embedding $\rho: |\CC|\to T$. Furthermore, convexity of the polyhedra $c\in Ob(\CC)$ implies that 
if $a, b\in \Ob(\CC)$ and $\rho_a(a)\subset \rho_b(\CP_b)$, then $a\le b$.

\begin{lemma}\label{proper}
1. All parasitic subspaces have dimension $\le N-2$. 

2. If $c$ is simple and the intersection 
$I_{c,\sigma}$ contains $\CP_{c'}$ for some face $c'$ of $\CP_{c}$, then 
$I_{c,\sigma}$ is not parasitic.  

3. Suppose that the complex $\CC$ is simple. Then no parasitic subspace $I_{c,\si,b}\subset \CP_b$ contains a face $a$ of $b$. 
\end{lemma}
\proof (1) is clear. 

(2) Without loss of generality we may assume that $\si=(c_1,...,c_k)$ is such that $c_1, c_2,...,c_k, c'$ are faces of $c$. 
Hence, $I_{c,\sigma}$ also contains $P_{c'}$ and for each $i=1,...,k$, $P_{c_i}$  contains the face $c'$. By convexity of $c$ it follows that 
$c'\subset c_i, i=1,...k$. Simplicity of $c$ then implies that for the face $c_{k+1}:= c_1\cap ...\cap c_k$ of $P_c$, we have 
$$
\Span(c_{k+1}) =\Span(c_1)\cap ...\cap \Span(c_k). 
$$ 
Thus, $I_{c,\sigma}= P_{c_{k+1}}$ and, hence, $I_{c,\sigma}$ is not parasitic.

(3) Let $I_{c,\si,b}=\rho_b^{-1} \rho_c (I_{c,\sigma})$ be parasitic in $\CP_b$. 
Assume that $a$ is a face of $b$ such that 
$a\subset I_{c,\si,b}$. Then $\rho_a(a)\subset \rho_c (I_{c,\sigma})\subset \rho_c(\CP_c)$ and, hence, $a\le c$. Now 
the assertion follows from (2). \qed

\medskip 
 We next  define a certain blow-up $\CB{\CP}$ of ${\CP}$ above which eliminates parasitic subspaces. 
In the process of the blow-up, the projective 
spaces $\CP_c\in \Ob(\im \CP)$ will be replaced with smooth rational varieties $B\CP_c$ so that the linear morphisms 
$F_{c_2,c_1}\in \mor(\im \CP)$ correspond to embeddings $bF_{c_2,c_1}: B\CP_{c_1}\to B\CP_{c_2}$.  The varieties $B{\mathcal P}$ are obtained by a sequence of  blow-ups of parasitic subspaces. 

\begin{prop}
Let $\CC$ be a simple Euclidean complex. Then there exists a strictly faithful complex of varieties $\CA: \CC\to \CB\CP(\CC)$ 
based on $\CC$ so that:

\begin{enumerate}
\item  The  direct limit of $im(\CA)=\CB\CP(\CC)$ exists and is a projective variety $X$ with simple normal crossing singularities.

\item $X$ is  an algebraic realization of $|\CC|$. 

\item $\pi_1(X)\cong \pi_1(|\CC|)$. 

\item If $\CC$ is rational, then $X$ is also defined over $\Q$.
\end{enumerate}
\end{prop}
\proof The complex $\CA$ is constructed by inductive blow-up of the complex $\CP$ above. 

We proceed by induction on dimension of parasitic intersections. First, for each face  
$c\in {\mathcal C}$ we blow up all parasitic subspaces in $\CP_c$ 
which are points. We will use the notation $B_0\CP_c$ for the resulting smooth 
rational varieties, $c\in \Faces(\CC)$. Observe that this blow-up is consistent with linear embeddings 
$F_{c_2,c_1}$, which, therefore, extend to embeddings $b_0F_{c_2,c_1}: B_0 \CP_{c_1}\to  B_0 \CP_{c_2}$. 
We let $\CB_0 {\mathcal P}$ denote the functor $\CC\to Varieties$, 
$$
c\mapsto B_0\CP_c, \quad f_{c_2,c_1}\mapsto b_0F_{c_2,c_1}. 
$$

Observe also that after the $0$-th blow-up, all 1-dimensional parasitic subspaces become pairwise disjoint (as we blew up their intersection points). We, thus, can now blow up each $B_0 \CP_c$ along every 1-dimensional blown-up parasitic subspaces   
$B_0 I_{c, \sigma}$. The result is a collection of smooth rational varieties $B_1 \CP_c$, $c\in \CC$. Again, the projective embeddings 
$b_0 F_{c_2,c_1}$ respect the blow-up, so we also get a collection of injective morphisms 
$b_1 F_{c_2,c_1}: B_1 \CP_{c_1}\to  B_1 \CP_{c_2}$. We, therefore, continue inductively on dimension of parasitic subspaces. 
After at most $N-1$ steps we obtain a complex $\CA=\CB {\mathcal P}$, whose image $\im(\CA)$ has blown-up projective spaces $B \CP_c$ as objects 
and embeddings $b F_{c_2, c_1}$ as morphisms. Observe, that the subvarieties along which we do the blow-up have dimension $\le N-2$. 
Moreover, we never have to blow up the entire $\CP_c$ for any $c\in Ob(\CC)$ (see Lemma \ref{proper}).  
Now, by the construction, the functor $\CA$ is strictly faithful. 

\medskip
It is easy to see that the conditions of Proposition \ref{seminorm.limit.thm} hold for $\CA$, the key is that the complex $\CC$ was simple 
and the normality condition was satisfied by the complex $\CP$. Therefore, the complex variety $X$ which is the direct limit 
of $im(\CA)$ exists and is an algebraic realization of $|\CC|$. We check in \ref{project.crit.say} that the variety $X$ is projective.

By the construction, since the complex $\CC$ is  simple, the variety $X$ has only normal crossing singularities. Since each $B\CP_c$ is 
simply-connected, Lemma \ref{pi_1-cong} implies that  $\pi_1(X)\cong \pi_1(|\CC|)$. Lastly, 
if we start with a rational complex $\CC$, then all the blow-ups are defined over $\q$ and so is the direct limit $X$.
\qed

\medskip
Suppose now that $W$ is a finite, connected, simplicial complex. By Corollary \ref{Hirsh+Cairns}, there exists a rational finite 
 simple polyhedral complex $\CC$ so that $|\CC|$ is homotopy-equivalent to $|W|$. Thus, we conclude

\begin{theorem}\label{thm:normcross}
There exists a complex projective variety $Z=Z_{\CC}$ defined over $\Q$ 
whose only singularities are simple normal crossings, so that $\pi_1(Z)\cong \pi_1(|W|)$.  
\end{theorem}

\section{Direct limits of complexes of varieties}\label{sec.5}

Example \ref{nolim.exmp} below shows that direct limits need not exist in
the category  of varieties
need not exist,
not even if all objects are smooth and all morphisms are closed embeddings.
By analyzing the example, we see that problems arise if
some of the images  $\phi_{ik}(X_i)\subset X_k$ and $\phi_{jk}(X_j)\subset X_k$
are tangent to each other but not if they are all transversal.
The right condition seems to be the seminormality of the images.

\begin{defn}\label{seminorm.defn}
 Recall that a complex space $X$ is called {\it normal}
if for every open subset $U\subset X$, every bounded meromorphic function
is holomorphic.

As a slight weakening, a complex space $X$ is called {\it seminormal}
if for every open subset $U\subset X$, every continuous meromorphic function
is holomorphic.

The following are some key examples:  
$(x^2=y^3)\subset \c^2$ and $(x^3=y^3)\subset \c^2$ are
 not   seminormal  (as shown by $x/y$ and  $x^2/y$) but
$(x^2=y^2)\subset \c^2$ is seminormal.

The key property that we use is the following.

Assume that $X$ is  seminormal.
Let $Y$ be any variety (or complex analytic space)
and $p:Y\to X$ be any algebraic  (or complex analytic morphism)  that is
a homeomorphism in the Euclidean topology. Then $p$ is an
isomorphism of varieties (or of complex analytic spaces).

\end{defn}

\begin{prop}\label{seminorm.limit.thm} Let 
 ${\mathbf X}:=\{ X_i: i\in I, \phi_{ij}:X_i\to X_j : (i,j)\in M\}$ 
 be a complex of varieties based on  a finite polyhedral complex $\CC$. 
 Assume that
 for each $k$ and each $J\subset I$
 the subvariety  $\cup_{j\in J}\im(\phi_{jk})\subset X_k$ is seminormal.
Then 
\begin{enumerate}
\item the direct limit ${\mathbf X}^{\infty}$ exists, 
\item the points of ${\mathbf X}^{\infty}$ are exactly the equivalence classes 
of points of $\amalg_{i\in I} X_i$, in particular, 
${\mathbf X}^{\infty}$ is an algebraic realization of $|\CC|$, and
\item  $\cup_{j\in J}\im(\phi_{j\infty})\subset {\mathbf X}^{\infty}$ is seminormal
for every $J\subset I$.

\end{enumerate}
\end{prop}

\proof The proof is by induction on $|I|$.

If there is a unique final object $X_j$ then  ${\mathbf X}^{\infty}=X_j$.

If not, let $X_j$ be a final object. Removing $X_j$ and all maps
to $X_j$, we get a smaller complex ${\mathbf Y}_j$.
  By induction it has a direct limit
  ${\mathbf Y}^{\infty}_j$. 

From ${\mathbf Y}_j$
 take away all the $X_k$ that do not map to $X_j$ and all maps
to such an $X_k$. Again we get a smaller complex ${\mathbf Z}_j$ 
whose  direct limit  is ${\mathbf Z}^{\infty}_j$.

There are maps  ${\mathbf Z}^{\infty}_j\to X_j$ and
 ${\mathbf Z}^{\infty}_j\to {\mathbf Y}^{\infty}_j$.
 We claim that these are both closed embeddings.

Both maps are clearly injective and their image is seminormal.
For the first this follows from our assumption  and
in the second case by induction and (3). 
As we noted in Definition \ref{seminorm.defn}, these imply that these maps 
 are  closed embeddings.

Now we claim that ${\mathbf X}^{\infty}$ is the universal push-out
$$
\begin{array}{ccc}
{\mathbf Z}^{\infty}_j&\to &X_j\\
\downarrow && \downarrow \\
{\mathbf Y}^{\infty}_j&\to & {\mathbf X}^{\infty}
\end{array}
$$
The existence of the push-out as an algebraic space
is proved in   \cite[Thm.3.1]{artin} 
and as a variety in \cite{ferrand},
see also  \cite[Cor.48]{k-quot}.
If a limit of a diagram of  seminormal varieties exists,
it is automatically seminormal. 

Finally we need to check that (3) holds. 
Let $W^{\infty}\subset X^{\infty}$ be the union of the images
of $\phi_{i\infty}(X_i)$ for $i\in I'$ for some $I'\subset I$.
Then  $W_j:=W^{\infty}\cap X_j$,  $W_j^Z:=W^{\infty}\cap {\mathbf Z}^{\infty}_j$ and
$W_j^Y:=W^{\infty}\cap {\mathbf Y}^{\infty}_j$ and
are all unions of images of some
of the $X_i$, hence these are  seminormal by induction.
Note that ${\mathbf Z}^{\infty}_j\cup W_j$ and
${\mathbf Z}^{\infty}_j\cup W^Y_j$ are also   unions of images of some
of the $X_i$, hence seminormal.

To check seminormality, we may assume that all varieties are affine.
Let $h$ be a continuous meromorphic function on $W^{\infty}$.
The restriction of $h$ to   $W_j^Z$ is again a 
continuous meromorphic function, hence it is holomorphic since
$W_j^Z$ is  seminormal. This restriction thus extends to a
 holomorphic  function $h^Z$ on ${\mathbf Z}^{\infty}_j$.

The restriction of $h$ to $W_j$ (resp.\ $W_j^Y$)
is a continuous meromorphic function, hence it is holomorphic since
$W_j^Z$ (resp.\ $W_j^Y$) is  seminormal.

Thus $h^Z$ and $h|_{W_j}$ define a  continuous meromorphic function
on ${\mathbf Z}^{\infty}_j\cup W_j$. It is thus  holomorphic 
and extends to a  holomorphic  function $h_j$ on $X_j$.
Similarly,  $h^Z$ and $h|_{W^Y_j}$  extend to a  holomorphic  function $h^Y_j$ on
${\mathbf Y}^{\infty}_j$.

Finally $h_j$ and  $h^Y_j$ agree on ${\mathbf Z}^{\infty}_j$,
thus, by the universality of the push-out, they
define a holomorphic function $h^{\infty}$ on $X^{\infty}$.
Its restriction to $W$ is $h$, hence $h$ is also holomorphic.\qed

\begin{say}[Projectivity]\label{project.crit.say}
 As the examples (\ref{pillows}) or \cite[Example 15]{k-quot} show,
even if the $X_i$ are all projective and the direct limit
${\mathbf X}^{\infty}$ exists, the latter need not be projective.
The main difficulty is the following.

Let $L^{\infty}$ be a line bundle on ${\mathbf X}^{\infty}$.
By pull-back we obtain line bundles $L_i$  and
isomorphisms $L_i\cong \phi_{ij}^*L_j$ with the expected
compatibility conditions.  $L^{\infty}$ is ample iff each
$L_i$ is ample. Conversely, giving   line bundles $L_i$  and
isomorphisms $L_i\cong \phi_{ij}^*L_j$ with the expected
compatibility conditions determines a line bundle $L^{\infty}$
on ${\mathbf X}^{\infty}$. The practical difficulty is that we need
to specify actual isomorphisms $L_i\cong \phi_{ij}^*L_j$; it is
not enough to assume that $L_i$ and $\phi_{ij}^*L_j$ are isomorphic.
For line bundles this is not natural to do since we usually specify
them only up-to the  fiberwise $\c^*$-action.

However, once we work only with  subsheaves of
 a fixed reference sheaf  $F^{\infty}$,
 the isomorphisms are easy to specify.

More generally,
 let $X=\cup_i X_i$ be a scheme with irreducible components
$X_i$. Let $F$ be a  coherent sheaf on $Y$. Then specifying a  coherent subsheaf
$G\subset F$ is equivalent to specifying   coherent subsheaves
$G_i\subset F|_{X_i}$ such that 
$G_i|_{X_i\cap X_j}=G_j|_{X_i\cap X_j}$ for every $i, j$. 
Furthermore, by the Nakayama lemma, $G$ is locally free iff
every $G_i$ is locally free.

In our case, we start with each irreducible component
identified with a $\p^n$ and then we blow up parasitic subvarieties.
Thus each   irreducible component $X_i$ comes with a natural morphism to $\p^n$.
For $F$ we choose the pull-back of $\o_{\p^n}(m)$ for some $m\gg 1$.
These pull-backs are not ample since they are trivial on the fibers
of $p_i:X_i\to \p^n$.

In Section \ref{sec.4} we construct the $X_i$ as follows.
\begin{enumerate}
\item Fix a smooth projective variety ${\mathbf P}$ with an 
ample line bundle $L$
(in our case in fact ${\mathbf P}\cong \p^n$ and $L\cong \o_{\p^n}(1)$).
Set $X_i^0:={\mathbf P}$.

\item If $X_i^j$ is already defined, we pick a smooth subvariety
$Z_i^j\subset X_i^j$ of dimension $j$ and let
$\pi_i^j:X_i^{j+1}\to X_i^j$ denote the blow-up of $Z_i^j$
with exceptional divisor $E_i^{j+1}$.

\item Set $X_i:=X_i^n$ with morphisms  $\Pi_i^j:X_i^n\to X_i^j$.
\end{enumerate}

{\it Claim. \ref{project.crit.say}.4.} For all  
$m_0\gg m_1\gg \cdots \gg m_n>0$,
the following line bundle is ample on $X_i$:
$$
\Bigl(\bigl(\Pi_i^0\bigr)^*H^{m_0}\Bigr)\Bigl(-m_1\bigl(\Pi_i^1\bigr)^*(E_i^1)-
\cdots - m_{n-1}\bigl(\Pi_i^1\bigr)^*(E_i^{n-1})-m_nE_i^n\Bigr)
$$

\proof 
Let $Y$ be a smooth variety and $Z\subset Y$ a smooth subvariety.
Let $p_Y:B_ZY\to Y$ denote the blow-up with exceptional divisor $E_Y$.
Let $H$ be an ample invertible sheaf on $X$.  
Then  $p_Y^*H^a( -b\cdot E_Y)$
is ample on $B_ZY$ for $a \gg b >0$ (cf. \cite[Prop.II.7.10]{hartsh}).

Applying this inductively to the blow-ups $\pi_i^j:X_i^{j+1}\to X_i^j$
we get our claim.\qed
\medskip

For later use, also not the following.
Assume that we have $Y_1\subset Y$ smooth and $Z_1:=Z\cap Y_1$
 also smooth.
Let $E_{Y_1}$ be the exceptional divisor of $p_{Y_1}:B_{Z_1}Y_1\to Y_1$.
Then there is an identity
$$
\bigl(p_Y^*H^a( -b\cdot E_Y)\bigr)|_{B_{Z_1}Y_1}=
p_{Y_1}^*\bigl(H|_{Y_1}\bigr)^a\bigl( -b\cdot E_{Y_1}\bigr).
$$

\medskip

(\ref{project.crit.say}.5) All the $X_i$ map to ${\mathbf P}$ in a
compatible manner, hence we have a fixed reference map
$\Pi^{\infty}: X^{\infty}\to {\mathbf P}$. For $m_0\gg 1$ we get 
our  reference sheaf  $F^{\infty}:=\bigl(\Pi^{\infty}\bigr)^* H^{m_0}$.
\medskip

(\ref{project.crit.say}.6) For each $i$ we have
$F_i:=F^{\infty}|_{X_i}=\bigl(\Pi_i^0\bigr)^* H^{m_0}$.
For fixed $m_0\gg m_1\gg \cdots \gg m_n>0$
the formula (\ref{project.crit.say}.4) defines a
subsheaf $G_i\subset F_i$ and $G_i$ is an ample line bundle
on $X_i$. As we noted above, all that remains is to prove that
$$
F^{\infty}|_{X_i\cap X_j }\supset G_i|_{X_i\cap X_j }=
G_j|_{X_i\cap X_j }\subset F^{\infty}|_{X_i\cap X_j }
\quad \forall i,j.
$$
This follows from the compatibility of blow-ups with
restrictions noted after the proof of (\ref{project.crit.say}.4).\qed
\end{say}
\medskip

The next example shows that direct limits need not exist in the
category of varieties. 

\begin{exmp}\label{nolim.exmp}
Start with the polyhedral subcomplex of $\r^2$ whose objects are
$$
(0,0), (x\leq 0, 0), (x\geq 0, 0),  (x, y\leq 0), (x, y\geq 0).
$$
We try to build  an algebraic realization with objects
$$
\c^1_x, \c^2_{x,y},  \c^2_{x,z}, \c^3_{\mathbf u}, \c^3_{\mathbf v}
\eqno{(\ref{nolim.exmp}.1)}
$$
where the maps are
$$
\begin{array}{l}
\c^1_x\to \c^2_{x,y}: x\mapsto (x,0)   \\
\c^1_x\to \c^2_{x,z}: x\mapsto (x,0)   \\
\c^2_{x,y}\to \c^3_{\mathbf u} : (x,y)\mapsto (x,y,0) \\
\c^2_{x,z}\to \c^3_{\mathbf u} : (x,z)\mapsto (x,z,z^2) \\
\c^2_{x,y}\to \c^3_{\mathbf v} : (x,y)\mapsto (x,y,0) \\
\c^2_{x,z}\to \c^3_{\mathbf v} : (x,z)\mapsto (x+z,z,z^2) 
\end{array}
\eqno{(\ref{nolim.exmp}.2)}
$$
These are all embeddings (even scheme theoretically).

We claim that if $X$ is any algebraic variety and
$g^i_*:\c^i_*\to X$ are algebraic maps from the $\c^i_*$ in (\ref{nolim.exmp}.1)
to $X$ with the expected compatibility properties
then $g^1_x:\c^1_x\to X$ is constant.
In an algebraic realization all the maps $g^i_*:\c^i_*\to X$
should be injective, thus this example is not an algebraic realization.
(A  careful analysis of the proof shows that the  direct limit does not exist
in the category of varieties, or even in the category of schemes of finite type.
The  direct limit  exists in the  category of  schemes
 but it is not Noetherian.)

So, take a regular function $\phi$ on $X$.
We can pull it back to $\c^2_{x,y}$ and $\c^2_{x,z}$ to get
 polynomials
$$
\tsum_{ij} a(i,j)x^iy^j\qtq{and} \tsum_{ij} b(i,j)x^iz^j.
$$
Next we compute these 2 ways.
First we pull $\phi$   back to
$\c^3_{\mathbf u}$. We get a polynomial 
$$
f(u_1, u_2, u_3)=\tsum_{ijk} c(i,j,k)u_1^i u_2^j u_3^k.
$$
Pull it back to $\c^2_{x,y}$ and $\c^2_{x,z}$  to get that
$$
a(i,j)= c(i,j,0)\qtq{and}  b(i,j)=c(i,j,0)+c(i,j-2,1)+c(i,j-4,2)+\cdots .
$$
Thus we obtain that
$$
a(i,0)=b(i,0)\qtq{and}  a(i,1)=b(i,1)\quad \forall\ i.
\eqno{(\ref{nolim.exmp}.3)}
$$
Next we pull $\phi$   back to
$\c^3_{\mathbf v}$. We get a polynomial 
$$
g(v_1, v_2, v_3)=\tsum_{ijk} d(i,j,k)v_1^i v_2^j v_3^k.
$$
Pull it back to $\c^2_{x,y}$  to get that
$a(i,j)= d(i,j,0)$. The pull-back to $\c^2_{x,z}$ 
involves the binomial coefficients; we are interested in the
first 2 terms only:
$$
 b(i,0)=d(i,0,0) \qtq{and} b(i,1)=d(i,1,0)+(i+1)d(i+1,0,0).
$$
Thus we obtain that
$$
a(i,0)=b(i,0)\qtq{and}  a(i,1)=b(i,1)-(i+1)b(i+1,0)\quad \forall\ i.
\eqno{(\ref{nolim.exmp}.4)}
$$
Comparing (\ref{nolim.exmp}.3) and (\ref{nolim.exmp}.4)
 we see that $a(i+1,0)=b(i+1,0)=0$ for $i\geq 0$,
that is, $\phi$ is constant on the image of $\c^1_x$.

Note that the same argument holds if $f,g$ are power series, thus
the problem is analytically local everywhere along $\c^1_x$.
In fact, the problem exists already
if we work with $C^2$-functions. (That is if $X\subset \r^N$ and we
require  $\c^3_{\mathbf u}\to X\subset \r^N$ and
$\c^3_{\mathbf v}\to X\subset \r^N$ to be at least $C^2$.)

\end{exmp}

\begin{exmp}[Triangular pillows]\label{pillows}
Take 2 copies of $\p^2_i:= \p^2(x_i:y_i: z_i)$ of $\c\p^2$ and 
the triangles  $C_i:=(x_iy_iz_i=0)\subset \p^2_i$.
Given $c_x, c_y, c_z\in \c^*$ 
define $\phi(c_x, c_y, c_z):C_1\to C_2$ by 
$(0:y_1:z_1)\mapsto (0:y_1:c_z z_1)$, $(x_1:0:z_1)\mapsto (c_xx_1:0: z_1)$
and $(x_1:y_1:0)\mapsto (x_1:c_yy_1: 0)$ 
and glue the 2 copies  of $\p^2$ 
using  $\phi(c_x, c_y, c_z)$ to get the surface
$S(c_x, c_y, c_z)$. 

We claim that $S(c_x, c_y, c_z)$
 is projective iff the product $c_xc_yc_z$ is a root of unity.

To see this  note that $\pic^0(C_i)\cong \c^*$
and  $\pic^r(C_i)$ is a principal homogeneous space
under $\c^*$ for every $r\in \z$. We can identify $\pic^3(C_i)$
 with $\c^*$ using
the restriction of the ample generator $L_i$ of 
$\pic\bigl(\p^2_i\bigr)\cong \z$ as the base point.

The key observation is that $\phi(c_x, c_y, c_z)^*:\pic^3(C_2)\to 
\pic^3(C_1)$ is the 
multiplication by $c_xc_yc_z$. 
Thus if $c_xc_yc_z$ is an $r$th root of unity then 
$L_1^r$ and $L_2^r$ glue together to an ample line bundle
 but otherwise every line bundle on $S(c_x, c_y, c_z)$
is topologically trivial.
\end{exmp}

\section{Proof of Theorem \ref{link.thm}}

So far, for every finitely presented group $G$ we have
constructed (Theorem \ref{thm:normcross}) a  complex projective variety $Z$ with  simple normal crossing singularities
  such that $\pi_1(Z)\cong G$. Using any such $Z$, we next construct
a singularity. This relies on the following result which is
a combination of 
\cite[Thm.8 and Prop.10]{kol-exmp}.

\begin{thm} \label{CI.sch.is.exc.set.thm}
Let $Z$ be an $(n\geq 2)$-dimensional 
 projective variety with  simple normal crossing singularities only
 and $L$ an ample line bundle on $Z$.
Then for $m\gg 1$ there are 
  germs of  normal  singularities  $\bigl(0\in X=X(Z,L,m)\bigr)$
with a partial resolution
$$
\begin{array}{ccc}
Z & \subset & Y\\
\downarrow && \hphantom{\pi}\downarrow\pi \\
0 & \in & X
\end{array}
\quad\qtq{where $Y\setminus Z\cong X\setminus\{0\}$}
$$
such that
\begin{enumerate}
\item  $Z$   is a Cartier divisor
in $Y$,
\item the normal bundle of $Z$ in  $Y$
is  $K_Z\otimes L^{-m}$,
\item $\pi_1\bigl(\res (X)\bigr)\cong \pi_1(Z)$, 
\item The kernel of $\pi_1\bigl(\link (X)\bigr)\onto \pi_1\bigl(\res (X)\bigr)$
is cyclic, central and generated by any loop around an irreducible component
of $Z$,
\item if  $\dim Z\leq 4$ then $(0\in X)$ is an isolated singular point. \qed
\end{enumerate}
\end{thm}

Note that the isomorphism  $\pi_1\bigl(\res (X)\bigr)\cong \pi_1(Z)$
is not explicitly stated in \cite[Thm.8]{kol-exmp}.
However, $Z$ is a deformation retract of $Y$, hence 
$\pi_1(Y)\cong \pi_1(Z)$.
By \cite[8.3]{kol-exmp} $Y$ has terminal singularities, hence,
by    \cite{k-shaf, tak},
$$
\pi_1\bigl(\res (X)\bigr)\cong\pi_1\bigl(\res (Y)\bigr)\cong \pi_1(Y)\cong 
\pi_1(Z).
$$
Alternatively, if $\dim Z=2$  (which is the only case that we need here),
we have a complete description of
the possible singularities of $Y$. By  \cite[Claim 5.10]{kol-exmp} they are
of the form  $(x_1x_2=x_3x_4)\subset \c^4$ with $x_3=0$ defining $Z$.
There are 2 local irreducible components of $Z$ given by
 $x_1=x_3=0$ and  $x_2=x_3=0$.
We can resolve these singularities by a single blow-up. 
The exceptional divisor is simply connected, hence
$\pi_1\bigl(\res (Y)\bigr)\cong \pi_1(Y)$.

Note also that under any such blow-up 
the dual simplicial complex changes by getting a new vertex on the edge
connecting the two local irreducible components.
Thus its homeomorphism type is unchanged.

It is clear that the kernel of 
$\pi_1\bigl(\link (X)\bigr)\onto \pi_1\bigl(\res (X)\bigr)$
is  generated by the loops around the irreducible components
of $Z$. Being cyclic, central is a special property of the actual construction
in \cite{kol-exmp}. However, our method  at the end of this Section 
 does not  use (\ref{CI.sch.is.exc.set.thm}.4).

In order to apply this to obtain an  isolated singular point,
we need a low dimensional variety $Z$. We use
 the following singular version of the
Lefschetz hyperplane theorem; see
\cite[Sec.II.1.2]{gm-book} for a stronger result and references.  

\begin{thm} Let $X$ be a projective variety of dimension $\geq 3$
with local complete intersection
 singularities and  $H\subset X$ a general hyperplane section.
Then $\pi_1(H)\cong \pi_1(X)$.
\end{thm}

Thus we can apply Theorem \ref{CI.sch.is.exc.set.thm}
with $\dim Z=2$. The remaining issue is
 to deal with the kernel
of $\pi_1\bigl(\link (X)\bigr)\onto \pi_1\bigl(\res (X)\bigr)$.
For this we make a more careful choice of $Z$.

Pick a smooth point $z\in Z$. We blow up $Z$ and then we blow up
a point on the exceptional curve to get $Z_2$. There are two exceptional
curves $E_1, E_2$ on $Z_2$ and  
$\bigl( K_{Z_2}\cdot E_1\bigr)=0$ and $\bigl( K_{Z_2}\cdot E_2\bigr)=-1$.
Since $\pi_1(Z_2)\cong \pi_1(Z)$, we can use this $Z_2$ in
 Theorem \ref{CI.sch.is.exc.set.thm} to get $Z_2\subset Y_2$.
Let $N$ be the normal bundle of $Z_2\subset Y_2$.

From (\ref{CI.sch.is.exc.set.thm}.2) and the above
computation we conclude that
$c_1(N)\cap [E_1]$ and $c_1(N)\cap [E_2]$ are relatively prime.

Since $Z_2$ and $Y_2$ are both smooth along $E_i$,
the boundary of the normal disc bundle restricted to $E_i$ is  a lens space 
$L_i$ with $|\pi_1(L_i)|=c_1(N)\cap [E_i]$.
Thus, if $\gamma_z$ denotes a small circle in $Y$ around $Z$
centered at $z$ then  the order of $\gamma_z$ in 
$\pi_1\bigl(\link (X)\bigr)$
divides both $c_1(N)\cap [E_1]$ and $c_1(N)\cap [E_2]$.
Hence the $\gamma_z$ is trivial in $\pi_1\bigl(\link (X)\bigr)$. 

We can perform these blow-ups for points in every irreducible component
of $Z$. Then every generator of the kernel of 
$\pi_1\bigl(\link (X)\bigr)\onto \pi_1\bigl(\res (X)\bigr)$
is trivial, hence $\pi_1\bigl(\link (X)\bigr)\cong \pi_1\bigl(\res (X)\bigr)$.
This completes the proof of
Theorem \ref{link.thm}. \qed

\section{Rational singularities and superperfect groups}
\label{sec.ratsing}

\begin{defn} \label{rtl.defn}
A quasi projective variety $X$ has {\it rational singularities}
if for one (equivalently every) resolution of singularities $p:Y\to X$
and for every algebraic (or holomorphic) vector bundle $F$ on $X$,
the natural maps  $H^i(X, F)\to H^i(Y, p^*F)$ are isomorphisms.
That is, for purposes of computing cohomology of vector bundles,
$X$ behaves like a smooth variety.
See \cite[Sec.5.1]{km-book} for details.
\end{defn}

\begin{defn}\label{DR.gen.defn}
Let $(0\in X)$ be a not-necessarily isolated  singularity and choose a
resolution of singularities $p:Y\to X$ such that $E:=p^{-1}(0)$
is a simple normal crossing divisor. 
Let $\{E_i\subset E: i\in I\}$ be the irreducible components
and for $J\subset I$ set $E_J:=\cap_{i\in J} E_i$.

The {\it dual simplicial complex} of $E$ has 
 vertices  $\{v_i:i\in I\}$ indexed by  the irreducible components  of $E$.
For $J\subset I$ we attach a $|J|$-simplex for 
every irreducible component of $\cap_{i\in J} E_i$.
Thus $D(E)$ is a  simplicial complex of dimension $\leq \dim X-1$.

The dual simplicial complex of a singularity seems to have been known to several
people but not explicitly studied until recently.
The dual graph of a normal surface singularity has a long history. Higher
dimensional versions appear in 
\cite{MR0506296, MR0466149, MR576865, friedman-etal}
but systematic investigations were started only recently; see
\cite{MR2320738, MR2399025, payne09, payne11}.

It is proved in \cite{ MR2320738, MR2399025} that  the homotopy type of $D(E)$
is  independent of the resolution $Y\to X$. As before, we denote it by
$\dres(0\in X)$. 

A possible argument runs as follows.
Let $F\subset E_J$ be an irreducible component.
If we blow up $F$, the dual simplicial complex
changes by a barycentric subdivision of the
$|J|$-simplex corresponding to $F$.
If we blow up a smooth subvariety $Z\subset F$ that is not contained
in any smaller $E_{J'}$ then 
 the dual simplicial complex
changes by attaching the cone over the star of  the
$|J|$-simplex corresponding to $F$.
Thus in both cases, the homotopy type of $D(E)$ is unchanged and
by \cite{wlod} this implies the general case.

If $X$ has rational singularities then 
$H^i\bigl(E, \o_E\bigr)=0$ for $i>0$ by \cite[2.14]{steenbrink}.
By Part 1 of Lemma \ref{friedman??.lem} below we conclude that 
$H^i\bigl(\dres(0\in X), \q\bigr)=0$ for $i>0$, That is,
$\dres(0\in X)$ is  $\q$-acyclic. 
\end{defn}

\begin{lem}\label{friedman??.lem}
Let $X$ be a simple normal crossing variety over $\c$ with
irreducible components $\{X_i:i\in I\}$. Let $T=D(X)$ be the dual 
simplicial complex of $X$. Then
\begin{enumerate}
\item There are natural injections
$H^r\bigl(T, \c\bigr)\into H^r\bigl(X, \o_{X}\bigr)$ for every $r$ and
\item 
For $J\subset I$ set $X_J:=\cap_{i\in J}X_i$ and
assume that   $H^r\bigl(X_J, \o_{X_J}\bigr)=0$  for every $r>0$ and
for every $J\subset I$.
Then  $H^r\bigl(X, \o_{X}\bigr)=H^r\bigl(T, \c\bigr)$ for every $r$.
\end{enumerate}
\end{lem}
\proof The following  is a combination of various arguments in
\cite[pp.68--72]{gri-sch} and
\cite[pp.26--27]{friedman-etal}.

Fix an ordering of $I$.
It is not hard to check that there is an exact  complex 
$$
0\to \c_X\to \tsum_i \c_{X_i}\to \tsum_{i<j} \c_{X_{ij}}\to \cdots
$$
where the $k$th term is
$ \tsum_{|J|=k} \c_{X_J}$ and $\c_{X_J}$ is the constant sheaf with support
$X_J$.   If $i\in J$ then the map
$\c_{X_{J\setminus i}}\to \c_{X_{J}}$ is the natural restriction
with a plus (resp.\ minus) 
sign if $i$ is in odd (resp.\ even) position in the ordering of $J$.

Thus the cohomology of $\c_X$ is also the hypercohomology of the
rest of the complex $\tsum_i \c_{X_i}\to \tsum_{i<j} \c_{X_{ij}}\to \cdots$.
This is computed by a spectral sequence whose $E_1$ term is
$$
\tsum_{|J|=q} H^p\bigl(X_J,  \c\bigr)\Rightarrow  H^{p+q}(X, \c).
\eqno{(\ref{friedman??.lem}.3)}
$$
The key observation is that this  spectral sequence
degenerates at $E_2$ \cite[p.??]{gri-sch}. The reason is that
$H^p\bigl(X_J,  \c\bigr) $ carries a Hodge structure of weight $p$
and there are no maps between  Hodge structures of weights.

Note also that the bottom (that is $p=0$)  row of (\ref{friedman??.lem}.3) is
$$
0\to  \tsum_i H^0\bigl(T_i,  \c\bigr)\to \tsum_{i<j}  H^0\bigl(T_{ij},  \c\bigr)
\to  \cdots
$$
where $T_i\subset T$ denotes the open star of the vertex corresponding
to $i\in I$ and  $T_J=\cap_{i\in J}T_i$.
The homology groups of this complex are exactly the $H^j\bigl(T,  \c\bigr)$.
Thus we have injections
$$
H^j\bigl(T,  \c\bigr)\into H^j\bigl(X,  \c_X\bigr).
\eqno{(\ref{friedman??.lem}.4)}
$$

Similarly, there is an exact  complex 
$$
0\to \o_X\to \tsum_i \o_{X_i}\to \tsum_{i<j} \o_{X_{ij}}\to \cdots
$$
which gives a spectral sequence whose $E_1$ term is
$$
\tsum_{|J|=q} H^p\bigl(X_J,  \o_{X_J}\bigr)\Rightarrow  H^{p+q}(X, \o_X).
\eqno{(\ref{friedman??.lem}.5)}
$$
By Hodge theory, the natural map from the  spectral sequence
(\ref{friedman??.lem}.5) to  the  spectral sequence
(\ref{friedman??.lem}.3) is a split surjection, hence
(\ref{friedman??.lem}.5) also degenerates at $E_2$ and so
$$
H^j\bigl(T,  \c\bigr)\into H^j\bigl(X,  \o_X\bigr).
\eqno{(\ref{friedman??.lem}.6)}
$$
is an injection. Under the assumptions of (2)
only the bottom row of (\ref{friedman??.lem}.5) is nonzero,
hence, in  this case,  
$H^j\bigl(T,  \c\bigr)= H^j\bigl(X,  \o_X\bigr)$. \qed

\medskip
In order to understand fundamental groups of links of rational singularities we need the following definition:

\begin{defn}\label{perfect}
Recall that a group $G$ is called {\em perfect} if 
it has trivial abelianization, equivalently, if $H_1(G,\z)=0$.
 Similarly, $G$ is called {\em superperfect} (see \cite{Berrick}) if $\tilde{H}_i(G,\Z)=0$ for $ i\le 2$. We generalize this notion to homology with coefficients in other commutative 
rings $R$: A group $G$ is $R$-perfect if  $H_1(G,R)=0$; $G$ is $R$-superperfect if $\tilde{H}_i(G,R)=0$ for $ i\le 2$. 
(We will be interested only in the cases $R=\Z$ and $R=\Q$.) 

Let $W$ be a cell complex.
Recall that by a theorem of  Hopf \cite{Hopf} the natural homomorphism
$H_2(W,R)\to H_2\bigl(\pi_1(W), R\bigr)$ is surjective and 
its kernel (in the case $R=\Z$) is the image of $\pi_2(W)$ under the Hurewicz homomorphism.

Therefore, if $\tilde{H}_i(W,R)=0$ for $ i\le 2$ then 
 $\tilde{H}_i\bigl(\pi_1(W),R\bigr)=0$ for $ i\le 2$.  
\end{defn}

To see surjectivity in Hopf's theorem observe the following: For $G=\pi_1(|W|)$ we let $f: W\to V=K(G,1)$ be the map inducing the isomorphism 
of fundamental groups. Then there exists a map of the 2-skeleta $h: V^{(2)}\to W^{(2)}$ which is a homotopy--right inverse to $f$. Hence, 
$H_2(f): H_2(W, R)\to H_2(V,R)=H_2(G,R)$ is onto for every commutative ring $R$.

\begin{exmp}\label{higman}
 Higman's group  $G=\<x_i| x_i [x_i, x_{i+1}], i\in \Z/4\Z\>$ 
is perfect, infinite and contains no proper finite index subgroups  \cite{Higman}. If $W$ is the (2-dimensional) 
presentation complex of $G$ then,  clearly, $\chi(W)=1$. Thus, $\tilde{H}_i(W,\Z)=0$, $i\le 2$.  In particular, 
$G$ is superperfect by Hopf's theorem. Moreover, $W$ is $K(G,1)$, see e.g. \cite{Bridson-Grunewald}.  
 Thus, $\tilde{H}_i(G,\Z)=0$ for all $i$.   
\end{exmp}

\begin{thm}\label{rtl.suprrperf.thm}
 Let $(0\in X)$ be a rational singularity. Then
$\pi_1\bigl(\dres(X)\bigr)$ is $\q$-superperfect and finitely presented. 
Conversely, for every finitely presented $\q$-superperfect group $G$
there is a 6-dimensional rational singularity $(0\in X)$ 
such that 
$$
\pi_1\bigl(\dres(X)\bigr)=\pi_1\bigl(\res(X)\bigr)=\pi_1\bigl(\link (X)\bigr)
\cong G.
$$
\end{thm}

\begin{rem}  
(1) The singularities constructed in Theorem \ref{rtl.suprrperf.thm}
are not isolated. Their singular locus is 1-dimensional.
Away from the origin 
it is the simplest possible non-isolated singularity,
locally  given by the equation
$$
\bigl(x_1^2+x_2^2+x_3^2+x_4^2+x_5^2+x_6^2=0\bigr)\subset \c^7.
$$
We do not know if in Theorem \ref{rtl.suprrperf.thm} one can get  
isolated singularities or not.


(2) For an arbitrary rational singularity $(0\in X)$,
the three groups $\pi_1\bigl(\dres(X)\bigr)$, $\pi_1\bigl(\res(X)\bigr)$
and $\pi_1\bigl(\link (X)\bigr) $ need not be isomorphic.
For example, if $\dim X=2$ then 
$ \pi_1\bigl(\dres(X)\bigr)=\pi_1\bigl(\res(X)\bigr)=1$
yet $\pi_1\bigl(\link (X)\bigr) $ can be infinite \cite{mumf}.

As another example, let $S$ be a fake projective plane,
that is, $H_i(S, \z)\cong H_i(\c\p^2, \z)$ for every $i$
yet $\pi_1(S)$ is infinite.
Such surfaces were classified in \cite{pra-yeu}.
Let $\bigl(0\in C(S)\bigr)$ denote a cone over $S$. Then
$ \pi_1\bigl(\dres(C(S))\bigr)=1$ yet
$\pi_1\bigl(\res(C(S))\bigr)=\pi_1(S)$ is infinite.
\end{rem}

{\em Proof of Theorem \ref{rtl.suprrperf.thm}.} The first claim of Theorem \ref{rtl.suprrperf.thm}
follows from the above cited results of \cite{steenbrink}
and  \cite{Hopf}. In order to see the converse,
for every finitely presented $\q$-superperfect group $G$ we
construct below (Theorem \ref{superperfect groups}) a simple 5-dimensional, $\q$-acyclic,  Euclidean polyhedral complex $\CC$ 
whose fundamental group is isomorphic to $G$.
Once this is done, we obtain a 5-dimensional projective variety $Z$ with
simple normal crossing singularities such that
$\pi_1(Z)\cong G$ and $H^i(Z, \o_Z)=0$ for $i>0$ by   Lemma \ref{pi_1-cong}.

We now apply  Theorem \ref{CI.sch.is.exc.set.thm}. 
The proof of  \cite[Prop.9.1]{kol-exmp} shows that for $m\gg 1$,
the resulting $X$ is a rational singularity.
As we noted, $X$ does not have isolated singularities
but they are completely described by 
\cite[Claim.5.10]{kol-exmp}. \qed

\medskip
Our next goal is to construct polyhedral complexes $\CC$ used in the proof of Theorem \ref{rtl.suprrperf.thm}.
The following theorem was proven by Kervaire for $R=\Z$, but examination of the proofs in \cite{Kervaire} and 
\cite{Kervaire-Milnor}  shows that they also apply to $R=\Q$.

\begin{theorem}\label{homology-sphere}
Every finitely-presented $R$-superperfect group is isomorphic to the fundamental group of a smooth $R$-homology 
$k$-sphere $M^k$ for every $k\ge 5$; here $R=\Z$ or $R=\Q$.  
\end{theorem}


\begin{corollary}\label{Z}
Let $R=\Z$ or $R=\Q$. Then a finitely-presented group $G$ is $R$-superperfect if and only if there exists a $5$-dimensional finite  
simplicial complex $W$ so that  $\pi_1(|W|)\cong G$ and $|W|$ is $R$-acyclic,
that is,  $\tilde{H}_*(W,R)=0$. 
\end{corollary}
\proof One direction of this corollary follows from Hopf's result above. Suppose that $G$ is $R$-superperfect and finitely-presented. 
Take the 5-dimensional homology sphere $M$ as in Theorem  
\ref{homology-sphere}. Since $M$ is smooth, we can assume that it is triangulated. 
Remove from $M$ the interior of a closed simplex. The result is the desired simplicial complex $W$. \qed 

\medskip
We now estimate the dimension of Euclidean thickening of $Z$ in Corollary \ref{Z}. 
A rough estimate is that $Z$ immerses in $\R^{10}$, since $Z$ is $5$-dimensional.  One can do much better as follows. 
Due to the results of \cite[\S 6]{Kervaire-Milnor}, the $5$-dimensional manifold $M^5$ constructed in Theorem \ref{homology-sphere} can be chosen to be {\em almost parallelizable}, i.e., the complement to a point $p$ in $M^5$ is parallelizable. Therefore, $M^5\setminus \{p\}$ admits an immersion in $\R^5$, see \cite{Phillips}. Hence, $W$ admits a 5-dimensional thickening $Y$, see 
Section \ref{thickening}. 
If $R=\Z$, then one can do even better and obtain a thickening $Y$ of $W$ which is an open subset of $\R^5$, see \cite{Livingston}.

In order to reduce the dimension of $X$ from $6$ to $5$ in Theorem \ref{rtl.suprrperf.thm} (and, thus, obtain isolated singularities) 
we have to impose further restrictions on the fundamental group $G$. 
Recall that a finite presentation of a group is called {\em balanced} if it has equal number of generators and relators. A group $G$ is called 
{\em balanced} 
if it admits a balanced presentation. Suppose that $G$ is an $R$-superperfect group which is the fundamental group of a 2-dimensional $R$-acyclic 
cell complex $W$. Without loss of generality, $W$ has exactly one vertex, i.e., $W$ is a presentation complex of $G$. 
Then $H_i(W,R)\cong H_i(G,R)=0, i=1,2$. In particular, $\chi(W)=1$. It then follows that $W$ has same number of 
edges and 2-cells. Hence, $G$ is balanced (with the balanced presentation complex $W$). 
Hausmann and Weinberger in \cite{HW} constructed examples of finite superperfect 
groups which are not balanced, see \cite{MR2118175}
 for more examples and a survey.
 Examples of finite $\Q$-superperfect groups which are not balanced are easier to construct: 
Take, for instance, the $k$-fold direct product $A_{p,k}=\Z/p\Z\times ... \times \Z/p\Z$  where $k\ge 2$. In particular, such groups do not admit 
$\Q$-acyclic presentation complexes
and they  do not occur as
$\pi_1\bigl(\dres(X)\bigr)$ for a 3-dimensional rational singularity.

Suppose that $G$ is balanced and $R$-superperfect ($R=\Z$ or $R=\Q$); then there exists a smooth $4$-dimensional $R$-homology sphere 
$M^4$ with the fundamental group $G$, see \cite{Kervaire}. Moreover, in Kervaire's construction one can assume that $M^4$ is almost parallelizable (i.e., it is a 4-dimensional Spin-manifold), see \cite{Kapovich}. Thus, for such $G$ there is a $4$-dimensional Euclidean thickening of its 2-dimensional (balanced) 
presentation complex $W$. More explicitly, in view of Stallings' theorem \cite{DR}, we can assume that $W$ embeds in $\R^4$. Since $W$ 
is a balanced  presentation complex of  a perfect group, $\chi(W)=0$, and, hence, $b_1(W)=b_2(W)=0$. 
Thus we obtain a $\Q$-acyclic 4-dimensional Euclidean thickening of $W$. 

Note that our methods cannot produce a 5-dimensional variety in Theorem \ref{rtl.suprrperf.thm} without the balancing condition. 
Specifically, given a $\Q$-superperfect group $G$ we would need a 4-dimensional $\Q$-acyclic manifold with the fundamental group $G$. 
However, one can show, repeating the arguments of \cite{HW}, that for all but finitely many finite groups $G$ constructed in \cite{HW} 
such 4-dimensional manifold does not exist. 

Reducing the thickening dimension to $3$ is, of course, very seldom possible since it amounts to assuming 
that $G$ is a 3-manifold group, which are quite rare among finitely-presented groups.

By combining these observations with Corollary \ref{Hirsh+Cairns}, we conclude 

\begin{theorem}\label{superperfect groups}
Let $G$ be an $R$-superperfect finitely-presented group ($R=\Z$ or $R=\Q$). Then there exists a finite  
simple $5$-dimensional Euclidean polyhedral complex $\CC$ so that $|\CC|$ is $R$-acyclic and has fundamental group 
isomorphic to $G$. Moreover, if $G$ admits a balanced presentation then we can take such $\CC$ to be $4$-dimensional. 
\end{theorem}

\begin{corollary}\label{cor:balanced}
Suppose that $G$ is a finitely-presented $\Q$-superperfect group which admits a balanced presentation. Then in 
Theorem \ref{rtl.suprrperf.thm} one can take $X$ which is $5$-dimensional and $(0\in X)$ an isolated singularity. 
\end{corollary}

\begin{corollary}\label{cor:trivial}
There exists  5-dimensional, isolated, rational singularities
 $(0\in X)$ so that the group 
$\pi_1^{\rm alg}(\link(X))$ is trivial yet  $\pi_1(\link(X))$ is infinite. 
\end{corollary}
\proof Take Higman's group $G$, see Example \ref{higman}. Then $G$ clearly has balanced presentation 
(its presentation has four generators and four relators), the group $G$ is also superperfect, infinite and 
has no nontrivial finite quotients. Now, the assertion follows from Theorem \ref{rtl.suprrperf.thm} and Corollary 
\ref{cor:balanced}. \qed

 \begin{ack} We are  grateful to Mike Freedman, Rob Kirby,
 Sam Payne and  Frank Quinn for references, to Kevin Walker for a sketch of the proof of 
 Corollary \ref{walker} and to Carlos Simpson for
many   comments and corrections. 
Partial financial support  to the first author was provided by the NSF grant number DMS-09-05802 and to 
the second author by  the NSF under grant number 
DMS-07-58275.
\end{ack}


\newcommand{\etalchar}[1]{$^{#1}$}
\providecommand{\bysame}{\leavevmode\hbox to3em{\hrulefill}\thinspace}
\providecommand{\MR}{\relax\ifhmode\unskip\space\fi MR }
\providecommand{\MRhref}[2]{%
  \href{http://www.ams.org/mathscinet-getitem?mr=#1}{#2}
}
\providecommand{\href}[2]{#2}

\medskip
Addresses:

\medskip 

\noindent M.K.: University of California, Davis, CA 95616 

{\begin{verbatim}kapovich@math.ucdavis.edu\end{verbatim}}

\medskip

\noindent J.K.: Princeton University, Princeton NJ 08544-1000

{\begin{verbatim}kollar@math.princeton.edu\end{verbatim}}

\end{document}